\documentstyle[12pt]{article}

\def\bl{\begin{eqnarray}}
\def\el{\end{eqnarray}}
\def\bll{\begin{eqnarray*}}
\def\ell{\end{eqnarray*}}
\def\su{\sum\limits}

 \title{On Relations among Fourier Coefficients and Sum-functions\thanks{ Supported in part by
Natural Science Foundation of China under grant number 10471130.}}
          \author{D. S. Yu and S. P. Zhou}

\date{}

\begin{document}
          \maketitle
          \pagenumbering{arabic}
          \begin{quote}
          \small
          \bf ABSTRACT: \rm We generalize five theorems of
          Leindler on the relations among  Fourier coefficients
          and sum-functions under the more general $NBV$
          condition.

                    \end{quote}

\begin{center}
\end{center}

\begin{flushleft}
{\bf \S 1. Introduction}
\end{flushleft}

Let $f(x)$ be a $p$ power integrable function of period $2\pi$, in symbol, $f\in L^p$, $p\geq 1$.
Define
$$\omega_p(f,h):=\sup\limits_{|t|\leq h}\|f(x+t)-f(x)\|_p$$
and
$$\omega_p^*(f,h):=\sup\limits_{0<t\leq h}\|f(x+t)+f(x-t)-2f(x)\|_p,$$
where $\|\cdot\|_p$ denotes the usual $L^p$ norm.

Denote by $E_n^{(p)}(f)$ the best approximation of order $n$ of $f$ in $L^p$.
The Lipschitz class $\Lambda_p$ and the Zygmund class
$\Lambda_{p}^*$ are defined by
$$\Lambda_p:=\{f\in L^p:\omega_p(f,h)=O(h)\}$$
and
$$\Lambda_p^*:=\{f\in L^p:\omega_p^*(f,h)=O(h)\}$$
respectively.

Leindler [5] introduced a class of sequences,  a natural extension of monotone decreasing sequences, named as $RBVS$.
Namely, a sequence ${\bf C}:=\{c_n\}$ of nonegative numbers
tending to zero is called of ``{\it rest bounded variation}", written as
${\bf C}\in RBVS$, if it satisfies
$$\sum\limits_{n=m}^\infty|\Delta c_n|\leq K({\bf C})c_m$$
for all $m=1,2,\cdots$, where $K({\bf C})$ is a constant only
depending upon ${\bf C}$ and $\Delta
c_n=c_n-c_{n+1},\;n=1,2,\cdots.$

Leindler [6] pointed out that $RBVS$ and the well known
quasi-monotone sequences ($CQMS$) are not comparable. Very recently, Le and
Zhou [2] defined a new condition named as $GBV$ condition
to include both the $RBV$ and quasi-monotone conditions. In the
special real case, the $GBV$ condition can be stated as follows: Let ${\bf
A}:=\{a_n\}$ be a sequence of nonegative numbers,
if
$$\sum\limits_{m=n}^{2n}|\Delta a_n|\leq K({\bf A})a_n$$
for all $n=1,2,\cdots,$ then we say ${\bf A}$ satisfies the $GBV$
condition, briefly, write ${\bf A}\in GBVS.$ Many important classic results
in Fourier analysis could be generalized  by replacing the monotonicity of
coefficients by $RBV$ or $GBV$ condition. For example, readers could refer to [11] for more information. Recently, we [8] further
introduced a new kind of sequences named as $NBVS$. In the real
case, the $NBVS$ can be defined as follows. Let ${\bf A}:=\{a_n\}$ be
a sequence of nonegative numbers, if
$$\sum\limits_{m=n}^{2n}|\Delta a_n|\leq K({\bf A})(a_n+a_{2n})$$
for all $n=1,2,\cdots,$ then we say ${\bf A}$ satisfies the $NBV$
condition, briefly, write ${\bf A}\in NBVS.$ As we know, the following embedding relations
$$RBVS\cup CQMS\subset GBVS\subset NBVS$$
holds.  Furthermore, as we mentioned in [8], $NBVS$ can be regarded as a ``two-sided"
monotonicity condition.

For convenience, through out the paper, we use $K$ to indicate a
positive constant which may depend upon $p$ and ${\bf A}$, its
value may be different even in the same line.

\begin{flushleft}
{\bf \S 2. Main Results }
\end{flushleft}

In this paper, we will establish the following results on the
relations among Fourier coefficients and the sum-functions. All of
them were proved for $RBVS$  by Leindler [7], and Theorem 1 was
proved for $GBVS$ by Zhou and Le [10].

\vspace{3mm} \bf Theorem 1. \quad\it Let ${\bf A}\in NBVS$ be
such for a fixed $p$, $1<p<\infty$, that
$$
\sum\limits_{n=1}^\infty n^{p-2}a_n^p<\infty.\hspace{.4in}(2.1)
$$
If $f$ is the sum of either of the series
$$
\sum\limits_{n=1}^\infty a_n\cos nx\;\;\;\mbox{or}\;\;\;
\sum\limits_{n=1}^\infty a_n\sin nx, \hspace{.4in}(2.2)$$
then
$$\omega_p(f,n^{-1})\leq
K_1n^{-1}\left\{\sum\limits_{\nu=1}^{n-1}\nu^{2p-2}a_\nu^p\right\}^{1/p}+K_2
\left\{\sum\limits_{\nu=n}^{\infty}\nu^{p-2}a_\nu^p\right\}^{1/p}.\hspace{.4in}(2.3)$$

\bf Theorem 2. \quad\it Let $1<p<\infty,\;1\leq
r<\infty$ and $\lambda(x)$, $x\geq 1$, be a positive monotone function
with $K_1\lambda(2^n)\leq\lambda(2^{n+1})\leq K_2\lambda(2^n)$,
where $K_1>0$. Write ${\bf A}=\{a_{n}\}\in NBVS$, $f(x)=\sum_{k=1}^{\infty}a_{k}\cos kx\in L^{p}$. Then
\bll
&&\sum\limits_{n=1}^\infty\lambda(n)a_n^r\leq KI(f,\lambda,r,p):=\nonumber\\
&&K\int_0^1\lambda\left(\frac{1}{t}\right)t^{r-2-\frac{r}{p}}
\left(\int_0^\pi|f(x+t)+f(x-t)-2f(x)|^pdx\right)^{r/p}dt.\hspace{.2in}(2.4)\ell
If $\sum_{n=1}^{\infty}|\Delta a_n|<\infty$ and $\lambda(x)$ satisfies the
additional conditions
$$\sum\limits_{n=1}^m\lambda(n)n^{\frac{r}{p}-r}\leq
K\lambda(m)m^{\frac{r}{p}-r+1}\hspace{.4in}(2.5)$$
and
$$
\sum\limits_{n=m}^\infty\lambda(n)n^{r(\frac{1}{p}-3)}\leq
K\lambda(m)m^{1+r(\frac{1}{p}-3)}, \hspace{.4in}(2.6)$$
then
$$
I(f,\lambda,r,p)\leq K\sum\limits_{n=1}^\infty\lambda(n)a_n^r.
\hspace{.4in}(2.7)$$

\bf Theorem 3. \quad\it Let $1<p<r$ and
$\{\varphi_n\}$ be a nonnegative nondecreasing sequence
satisfying
$\varphi_{n^2}\leq
K\varphi_n$ for all $n$. Define
$$\Phi(x):=\sum\limits_{n=1}^xn^{\frac{r}{p}-2}\varphi_n,$$
where $\varphi(x):=\varphi_n$ if $x\in (n-1,n).$
Write ${\bf A}=\{a_{n}\}\in NBVS$, $f(x)=\sum_{k=1}^{\infty}a_{k}\cos kx\in L^{p}$. Then the
statements
$$
\sum\limits_{n=1}^\infty \varphi_nn^{r-2}a_n^r<\infty,$$
$$
\sum\limits_{n=1}^\infty
\varphi_nn^{rs+\frac{r}{p}-2}\left(\sum\limits_{k=1}^nk^{(s+1)p-2}a_k^p\right)^{r/p}<\infty
\;\; \mbox{for any } \; s>\frac{1}{p}-\frac{1}{r},$$
$$\sum\limits_{n=1}^\infty
\varphi_nn^{\frac{r}{p}-2}\left(\sum\limits_{k=n}^\infty
k^{p-2}a_k^p\right)^{r/p}<\infty,$$
$$
\sum\limits_{n=1}^\infty
\varphi_nn^{\frac{r}{p}-2}\left(\omega_p\left(f,\frac{1}{n}\right)\right)^r<\infty,$$
$$\sum\limits_{n=1}^\infty
\varphi_nn^{\frac{r}{p}-2}\left(E_n^{(p)}(f)\right)^r<\infty,$$
$$
\int_0^\pi|f(x)|^{r-\frac{r}{p}+1}\Phi(|f(x)|)dx<\infty,
$$
$$\int_0^\pi|f(x)|^{r}\varphi(|f(x)|)dx<\infty,
$$
$$\int_0^\pi|f(x)|^{r}\varphi\left(\frac{1}{x}\right)dx<\infty
$$
and
$$
\int_0^\pi\varphi\left(\frac{1}{t}\right)t^{-\frac{r}{p}}
\left(\int_0^\pi|f(x+t)+f(x-t)-2f(x)|^pdx\right)^{r/p}<\infty
$$are equivalent.

\vspace{3mm} \bf Theorem 4. \quad\it If $1<p<\infty$ and ${\bf
A}\in NBVS$, then
$$
\sum\limits_{n=1}^\infty n^{2p-2}a_n^p<\infty\hspace{.4in}(2.8)$$
is a necessary and sufficient condition that a sum-function of
either of the series $(2.2)$

$(i)$ belongs to $\Lambda_p$, or

$(ii)$ is equivalent to an absolutely continuous function whose
derivative belongs to $L^p$.

\vspace{3mm} \bf Theorem 5. \quad\it If $f\in L^p$, $1<p<\infty$,
${\bf A}\in NBVS$ and $f$ is a sum-function of either the
series $(2.2)$, then $f\in\Lambda_p^*$ implies that
$$
\omega_p(f,h)\leq Kh|\log h|^{1/p}.
\hspace{.4in}(2.9)$$

\rm The proofs of the above results will be proceeded as in a way as
those of Leindler [7], only necessary modifications will be noted.

\begin{flushleft}
{\bf \S 3. Lemmas }
\end{flushleft}

\vspace{3mm} \bf Lemma 1. \quad\it Let $1<p<\infty,\;1\leq
r<\infty$ and $\lambda(x)$, $x\geq 1$, be a positive monotone function
with $K_1\lambda(2^n)\leq\lambda(2^{n+1})\leq K_2\lambda(2^n)$,
where $K_1>0$. Let $a_{n}\geq 0$, $n=1, 2, \cdots$, and write $f(x)=\sum_{k=1}^{\infty}a_{k}\cos kx\in L^{p}$. Then
$$
\sum\limits_{n=1}^\infty\lambda(n)n^{-r}\left(\sum\limits_{k=[n/2]}^{2n}a_k\right)^r\leq
KI(f,\lambda,r,p).\hspace{.4in}(3.1)$$
If $\sum\limits_{n=1}^\infty|\Delta a_n|<\infty$ and $\lambda(x)$
satisfies $(2.5)$ and $(2.6)$, then
$$I(f,\lambda,r,p)\leq K\sum\limits_{n=1}^\infty\lambda(n)
\left(\sum\limits_{k=n}^\infty|\Delta a_k|\right)^r. \hspace{.4in}(3.2)$$

{\bf Proof.} \rm The second result can be found directly in [4],
while the first can be proved by the same argument as of [4].

 \vspace{3mm} \bf Lemma 2([9]). \quad\it Let
$1<p<\infty$, $\{a_n\}\in NBVS$, and $f$ be the sum of either of the series $(2.2)$, then $f\in L^p$ if and only if
$(2.1)$ holds.

\vspace{3mm} \bf Lemma 3([3]). \quad\it Let $\alpha_n\geq 0$ and
$\lambda_n\geq 0$ be given, $\nu_1<\cdots<\nu_n<\cdots$ denote
the indices for which $\lambda_{\nu_n}>0$, and $N$ denote the
number of positive terms of the sequence $\lambda_n$, provides
this number is finite, or in the contrary case set $N=\infty.$ Set
$\nu_0=0$ and if $N<\infty$ then $\nu_{N+1}=\infty.$ We have the
following inequalities:
$$\sum\limits_{n=1}^\infty\lambda_n\left(\sum\limits_{k=1}^n
\alpha_k\right)^p\leq
p^p\sum\limits_{n=1}^\infty\lambda_{\nu_n}^{1-p}\left(\sum\limits_{k=\nu_n}^\infty
\lambda_k\right)^p\left(\sum\limits_{k=\nu_{n-1}+1}^{\nu_n}\alpha_k\right)^p,\hspace{.4in}(3.3)$$
$$
\sum\limits_{n=1}^\infty\lambda_n\left(\sum\limits_{k=n}^\infty
\alpha_k\right)^p\leq
p^p\sum\limits_{n=1}^N\lambda_{\nu_n}^{1-p}\left(\sum\limits_{k=1}^{\nu_n}
\lambda_k\right)^p\left(\sum\limits_{k=\nu_{n}}^{\nu_{n+1}-1}\alpha_k\right)^p. \hspace{.4in}(3.4)$$

 \bf Lemma 4([9]). \quad\it If ${\bf A}:=\{a_n\}\in
NBVS$, then for all $n\geq 1$, it holds that
$$I:=\sum\limits_{k=n}^\infty|\Delta a_k|\leq C({\bf A})
\left(a_n+a_{2n}+a_{4n}+\sum\limits_{k=n}^\infty
\frac{a_k}{k}\right).$$

 \bf Lemma 5. \quad\it If ${\bf A}:=\{a_n\}\in NBVS$,
then $$
n^{-p}\sum\limits_{m=1}^{n-1}m^{-2}\left(\sum\limits_{\nu=1}^m\nu^2|\Delta
a_\nu|\right)^p\leq
K\left(n^{-p}\sum\limits_{\nu=1}^{n-1}\nu^{2p-2}a_\nu^p+
\sum\limits_{\nu=n}^{\infty}\nu^{p-2}a_\nu^p\right).$$

\bf Proof. \rm Write
$$N(x)=:\frac{\log x}{\log 2},\;\;\;x>0.$$
Since
$m/4\leq 2^{[N(m/2)]}\leq m/2$ for $m\geq 2$,
then
\bll\su_{\nu=1}^m\nu^2|\Delta a_\nu|\leq
\su_{k=1}^{[N(m/2)]}\su_{\nu=2^{k-1}}^{2^k}\nu^2|\Delta
a_\nu|+\su_{k=[m/4]+1}^m\nu^2|\Delta a_\nu|=:J_1+J_2.\ell
 By the definition of
$NBVS$, we get
$$a_{2^k}\leq \su_{i=s}^{2^k}|\Delta a_i|+a_s\leq \su_{i=s}^{2s}|\Delta a_i|+a_s
\leq K(a_s+a_{2s}),$$
 and
$$a_{2^{k-1}}\leq \su^{s-1}_{i=2^{k-1}}|\Delta a_i|+a_s\leq K(a_{[s/2]}+a_{s})$$
for all $2^{k-1}\leq s\leq 2^k$,  and hence deduce that \bll J_1&\leq&
\su_{k=1}^{[N(m/2)]}2^{2k}\su_{\nu=2^{k-1}}^{2^k}|\Delta
a_\nu|\leq
K\su_{k=1}^{[N(m/2)]}2^{2k}\left(a_{2^{k-1}}+a_{2^k}\right)\\
&\leq&K\su_{k=1}^{[N(m/2)]}2^{k}\su_{s=2^{k-1}}^{2^k}\left(a_{[s/2]}+a_s+a_{2s}\right)
\leq K\su_{k=1}^{m}ka_k. \ell

 For any
 $[m/4]+1\leq s\leq 2[m/4]+2$,
 if $s$
is an even number, it follows from the definition of $NBVS$ that
\bll |a_{[m/4]+1}|&\leq&\su_{k=[m/4]+1}^{s-1}|\Delta
a_k|+a_s\leq \su_{k=[s/2]}^{s}|\Delta a_k|+a_s\\
&\leq&K(a_s+a_{[s/2]}). \ell
 If $s$ is an odd number, then \bll
|a_{[m/4]+1}|&\leq&\su_{k=[m/4]+1}^{s-1}|\Delta
a_k|+a_s\leq \su_{k=[s/2]}^{s-1}|\Delta a_k|+a_s\\
&\leq&K(a_s+a_{s-1}+a_{[s/2]}).\ell Therefore, in any case,
 $$m^2a_{[m/4]+1}\leq K\su_{s=[m/4]+1}^{2[m/4]+2}s(a_s+a_{s-1}+a_{[s/2]})\leq
K\su_{k=1}^{m}ka_k.\hspace{.4in}(3.5)$$
A similar discussion leads to
$$m^2|a_{2[m/4]+2}|\leq
K\su_{k=[m/4]+2}^{4[m/4]}s(a_s+a_{s-1}+a_{[s/2]})\leq
K\su_{k=1}^{m}ka_k, \hspace{.4in}(3.6)$$
 and
$$m^2|a_{[m/2]}|\leq
K\su_{k=[m/2]}^{m}s(a_s+a_{s-1}+a_{[s/2]})\leq
K\su_{k=1}^{m}ka_k. \hspace{.4in}(3.7)$$

Set
$$m^*=:\left\{\begin{array}{ll}m,&m\;\mbox{is even,}\\
m-1,&m\;\mbox{is odd.}\end{array}\right.$$ By (3.5)-(3.7), we
deduce that \bll J_2&\leq&\su_{\nu=[m/4]+1}^{2[m/4]+2}\nu^2|\Delta
a_\nu|+\su_{\nu=[m/2]}^{m^*}\nu^2|\Delta
a_\nu|+m^2(a_m+a_{m+1})\\
&\leq&Km^2\left(a_{[m/4]+1}+a_{2[m/4]+2}+a_{m^*/2}+a_{m^*}+a_m+a_{m+1}\right)\\
&\leq&Km^2\left(a_{[m/4]+1}+a_{2[m/4]+2}+a_{[m/2]}+a_{m-1}+a_m+a_{m+1}\right)\\
&\leq&K\su_{k=1}^{m}ka_k+Km^2\left(a_{m-1}+a_m+a_{m+1}\right).\ell
Combining all the estimates for $J_1$ and $J_2$ with the fact (see [7])
$$n^{-p}\sum\limits_{m=1}^{n-1}m^{-2}\left(\sum\limits_{\nu=1}^m\nu
a_\nu\right)^p\leq
Kn^{-p}\sum\limits_{\nu=1}^{n-1}\nu^{2p-2}a_\nu^p,$$
 we see
that \bll
n^{-p}\sum\limits_{m=1}^{n-1}m^{-2}\left(\sum\limits_{\nu=1}^m\nu^2|\Delta
a_\nu|\right)^p\leq
Kn^{-p}\sum\limits_{m=1}^{n-1}m^{-2}\left(\su_{k=1}^{m}ka_k\right)^p\ell
\bll\hspace{1cm}+Kn^{-p}\sum\limits_{m=1}^{n-1}m^{2p-2}(a_{m-1}^p+a_m^p+a_{m+1}^p)\ell
\bll&\leq& Kn^{-p}\sum\limits_{\nu=1}^{n-1}\nu^{2p-2}a_\nu^p+
Kn^{-p}\sum\limits_{m=1}^{n-1}m^{2p-2}a_{m+1}^p\\
&\leq& Kn^{-p}\sum\limits_{\nu=1}^{n-1}\nu^{2p-2}a_\nu^p+
Kn^{p-2}a_n^p\\
&\leq& Kn^{-p}\sum\limits_{\nu=1}^{n-1}\nu^{2p-2}a_\nu^p+K
\sum\limits_{\nu=n}^{\infty}\nu^{p-2}a_\nu^p.\ell

 \bf Lemma 6. \quad\it If ${\bf A}:=\{a_n\}\in NBVS$,
then
 \bll
n^{-p}\sum\limits_{m=1}^{n-1}m^{p-2}\left(\sum\limits_{\nu=m+1}^n\nu|\Delta
a_\nu|\right)^p\leq
K\left(n^{-p}\sum\limits_{\nu=1}^{n-1}\nu^{2p-2}a_\nu^p+
\sum\limits_{\nu=n}^{\infty}\nu^{p-2}a_\nu^p\right).\nonumber\\\ell

{\bf Proof.} \rm First assume that $n\geq 8(m+1)$. By noting that $m+1\leq
2^{[N(m+1)]+1}\leq 2(m+1)$ and $n/4\leq 2^{[N(n/2)]}\leq n/2$, we
can split $\su_{k=m+1}^n\nu|\Delta a_\nu|$ into
$$\su_{\nu=m+1}^n\nu|\Delta a_\nu|\leq\su_{k=[N(m+1)]+3}^{N(n/2)}\su_{\nu=2^{k-1}}^{2^k}
\nu|\Delta a_\nu|+\su_{\nu=[n/4]+1}^n\nu|\Delta
a_\nu|+\su_{\nu=m+1}^{2^{[N(m+1)]+3}-1}\nu|\Delta a_\nu|$$
$$\hspace{2cm}=:H_1+H_2+H_3.$$
Similar to what we have done for $J_1$ in the proof of Lemma 5, we get
$$H_1\leq K\su_{\nu=m+1}^n a_\nu.$$
Since $n\geq 8(m+1)$, then setting
$$N^{*}=:\left[\frac{1}{2}([n/4]+1)\right]\geq m+1,$$
again similar to $J_2$, we have
$$H_2\leq
K\su_{k=N^*}^{n}a_k+Kn\left(a_{n-1}+a_n+a_{n+1}\right)$$
$$ \leq
K\su_{k=m+1}^{n}a_k+Kn\left(a_{n-1}+a_n+a_{n+1}\right).$$ At the same time,
it is evident that
$$H_3\leq \su_{\nu=m+1}^{8(m+1)}\nu|\Delta a_\nu|\leq Km(a_{m+1}+a_{2(m+1)}+a_{4(m+1)}
+a_{8(m+1)}).$$ In case $n<8(m+1)$, then it clearly holds that
$$\su_{\nu=m+1}^n\nu|\Delta a_\nu|\leq Km(a_{m+1}+a_{2(m+1)}+a_{4(m+1)}+a_{8(m+1)}).$$
Altogether, all the above estimates lead to that \bll
\lefteqn{n^{-p}\sum\limits_{m=1}^{n-1}m^{p-2}\left(\sum\limits_{\nu=m+1}^n\nu|\Delta
a_\nu|\right)^p\leq
Kn^{-p}\sum\limits_{m=1}^{n-1}m^{p-2}\left(\su_{\nu=m+1}^n
a_\nu\right)^p}\\
&&+Kn^{p-1}(a_{n-1}^p+a_{n}^{p}+a_{n+1}^{p})+Kn^{-p}\sum\limits_{m=1}^{n-1}m^{2p-2}
\su_{j=1}^8a^{p}_{j(m+1)}.\ell
Obviously, $$a_{n}\leq \su_{i=k}^{n-1}|\Delta a_i|+a_k\leq K(a_k+a_{2k})$$
for $[n/2]+1\leq k\leq n$, which implies that
$$a_n\leq
Kn^{-1}\su_{k=[n/2]+1}^{n-1}(a_k+a_{2k})\leq
Kn^{-1}\su_{k=[n/2]+1}^{2n-2}a_k,\hspace{.4in}(3.8)$$
so that applying H\"{o}lder's inequality yields that
$$n^{p-1}a_{n}^{p}\leq Kn^{-1}
\left(\sum_{k=[n/2]+1}^{2n-2}a_{k}\right)^{p}
\leq Kn^{p-2}\sum_{k=[n/2]+1}^{2n-2}a_{k}^{p}$$
$$\leq K\left(n^{-p}\sum_{k=[n/2]+1}^{n-1}k^{2p-2}a_{k}^{p}+
\sum_{k=n}^{2n-2}k^{p-2}a_{k}^{p}\right)\leq K\left(n^{-p}
\sum_{k=1}^{n-1}k^{2p-2}a_{k}^{p}+\sum_{k=n}^{\infty}k^{p-2}a_{k}^{p}\right).
\hspace{.2in}(3.9)$$
Similarly,
$$n^{p-1}a_{n-1}^{p}\leq K\left(n^{-p}\sum_{k=1}^{n-1}k^{2p-2}a_{k}^{p}+\sum_{k=n}^{\infty}k^{p-2}a_{k}^{p}\right),\hspace{.2in}(3.10)$$
$$n^{p-1}a_{n+1}^{p}\leq K\left(n^{-p}\sum_{k=1}^{n-1}k^{2p-2}a_{k}^{p}+\sum_{k=n}^{\infty}k^{p-2}a_{k}^{p}\right).\hspace{.2in}(3.11)$$
Since (see [7]) \bll
n^{-p}\sum\limits_{m=1}^{n-1}m^{p-2}\left(\su_{\nu=m+1}^n
a_\nu\right)^p\leq Kn^{-p}\su_{\nu=1}^{n-1}\nu^{2p-2}a_\nu^p,\ell
with the estimates (3.9)-(3.11), Lemma 6 will be completed  if  we can verify that $$
\Delta=:n^{-p}\sum\limits_{m=1}^{n-1}m^{2p-2}
\su_{j=1}^8 a_{j(m+1)}^p\leq
Kn^{-p}\su_{\nu=1}^{n-1}\nu^{2p-2}a_\nu^p+\su_{\nu=n}^\infty\nu^{p-2}a_\nu^p.\hspace{.4in}(3.12)$$
Indeed, we prove (3.12) by the following way:
 \bll
\Delta&=&n^{-p}\su_{j=1}^8\su_{\nu=1}^{[n/j]-1}\nu^{2p-2}a_{j(\nu+1)}^p+
n^{-p}\su_{j=1}^8\su_{\nu=[n/j]}^{n-1}\nu^{2p-2}a_{j(\nu+1)}^p\\
&\leq&K\left(n^{-p}\su_{\nu=1}^{n-1}\nu^{2p-2}a_\nu^p
+n^{p-2}a_n^p+\su_{j=1}^8\su_{\nu=n}^{jn}\nu^{p-2}a_\nu^p\right)\\
&\leq&K\left(n^{-p}\su_{\nu=1}^{n-1}\nu^{2p-2}a_\nu^p
+\su_{j=1}^8\su_{\nu=n}^{jn}\nu^{p-2}a_\nu^p\right)\\
&\leq&K\left(n^{-p}\su_{\nu=1}^{n-1}\nu^{2p-2}a_\nu^p+
\su_{\nu=n}^\infty\nu^{p-2}a_\nu^p\right).\ell

\begin{flushleft}
{\bf \S 4. Proofs }
\end{flushleft}

{\bf Proof of Theorem 1.} \rm By Lemma 2, we know that the
condition (2.1) is both necessary and sufficient for
$f\in L^p$, $p>1$. We only need to treat the cosine series case, the other could be done similarly. Assume that $h=\pi/2n$. By the
symmetry of $f$, it is clear that
\begin{eqnarray*}
\omega_p(f,h)\leq K\sup\limits_{0<t\leq
h}\Big(\left\{\int_0^{\pi/n}|f(x\pm t)-f(x)|^pdx
\right\}^{1/p}\\
+\left\{\int_{\pi/n}^\pi|f(x\pm t)-f(x)|^pdx
\right\}^{1/p}\Big):=K\sup\limits_{0<t\leq h}(I_1+I_2).
\end{eqnarray*}
As the way done by Leindler [7], we have
$$
\frac{1}{2}I_1
\leq t\left\{\int_0^{\pi/n}\left(\sum\limits_{\nu=1}^{n-1}\nu
a_\nu\right)^pdx\right\}^{1/p}+K\left\{\sum\limits_{m=n}^\infty\int_{3\pi/2(m+1)}^{3\pi/2m}
\left|\sum\limits_{\nu=n}^\infty a_\nu\cos\nu
x\right|^pdx\right\}^{1/p}$$
$$\hspace{1.0cm}=:I_{11}+I_{12}.$$
Applying H\"{o}lder's inequality leads to
$$I_{11}\leq Kn^{-1}\left\{\sum\limits_{\nu=1}^{n-1}\nu^{2p-2}a_\nu^p\right\}^{1/p}.$$
By Abel's transformation and Lemma 4, we have
\begin{eqnarray*}\left|\sum\limits_{\nu=n}^\infty a_\nu\cos\nu
x\right|&\leq&
\sum\limits_{\nu=n}^m
a_\nu+(m+1)\sum\limits_{\nu=m+1}^\infty|\Delta a_\nu|\\
&\leq& K\left(\sum\limits_{\nu=n}^m
a_\nu+m(a_m+a_{2m}+a_{4m})+m\sum\limits_{\nu=m+1}^{\infty}\frac{a_\nu}{\nu}\right).
\end{eqnarray*}
Setting $\lambda_m=m^{-2}$ and $\alpha_m=0$ for $m<n$ and
$\alpha_m=a_m$ for $m\geq n$, we get
\begin{eqnarray*}
\sum\limits_{m=n}^\infty
m^{-2}\left(\sum\limits_{\nu=n}^m\alpha_\nu\right)^p&=&
\sum\limits_{m=1}^\infty m^{-2}
\left(\sum\limits_{\nu=1}^m\alpha_\nu\right)^p\\
&\leq& K\sum\limits_{m=1}^\infty
m^{p-2}\alpha_m^p=K\sum\limits_{m=n}^\infty
m^{p-2}a_m^p\end{eqnarray*} by (3.3). Again setting
$$\nu_1=n,\;\nu_2=n+1,\cdots,\nu_j=n+j,\cdots,$$
$$\lambda_1=\lambda_2=\cdots=\lambda_{\nu_1-1}=0,\;\lambda_{\nu_j}=\nu_j^{p-2},
\;j=1,2,\cdots,$$ with (3.4), we get
\begin{eqnarray*}
\su_{m=n}^\infty
m^{p-2}\left(\su_{\nu=m}^\infty\frac{a_\nu}{\nu}\right)^p&=&
\sum\limits_{j=1}^\infty\lambda_j\left(\sum\limits_{k=j}^\infty\frac{a_k}{k}\right)^p\nonumber\\
&\leq&
p^p\sum\limits_{j=1}^\infty\lambda_{\nu_j}^{1-p}\left(\sum\limits_{k=1}^{\nu_j}\lambda_k
\right)^p\left(\sum\limits_{k=\nu_j}^{\nu_{j+1}-1}\frac{a_k}{k}\right)^p\nonumber\\
&=&p^p\sum\limits_{m=n}^{\infty}m^{(p-2)(1-p)}\left(\sum\limits_{k=n}^mk^{p-2}\right)^p\left(
\frac{a_m}{m}\right)^p\nonumber\\
&\leq&p^p\sum\limits_{m=n}^\infty m^{p-2}a_m^p.\hspace{.4in}(4.1)
\end{eqnarray*}
 Therefore
\begin{eqnarray*}
I_{12}^p&\leq& K\sum\limits_{m=n}^\infty
m^{-2}\left(\sum\limits_{\nu=n}^m
a_\nu\right)^p+K\sum\limits_{m=n}^\infty
m^{p-2}a_m^p+K\sum\limits_{m=n}^\infty
m^{p-2}\left(\sum\limits_{\nu=m}^\infty \frac{a_\nu}{\nu}\right)^p\\
&\leq& K\sum\limits_{m=n}^\infty m^{p-2}a_m^p.
\end{eqnarray*}
Let $D_\nu(x)$ be the Dirichlet Kernel. Following the way of
Leindler [7], we see that
\begin{eqnarray*}
I_2\leq \left\{\int_{\pi/n}^\pi\left|\sum\limits_{\nu=1}^n\Delta
a_\nu[D_\nu(x\pm
t)-D_\nu(x)]\right|^pdx\right\}^{1/p}\\
+\left\{\int_{\pi/n}^\pi\left|\sum\limits_{\nu=n+1}^\infty\Delta
a_\nu[D_\nu(x\pm
t)-D_\nu(x)]\right|^pdx\right\}^{1/p}:=I_{21}+I_{22},
\end{eqnarray*}
and
\begin{eqnarray*}
I_{21}^p\leq
K\sum\limits_{m=1}^{n-1}\int_{\pi/(m+1)}^{\pi/m}\sum\limits_{\nu=1}^{n}|\Delta
a_\nu[D_\nu(x\pm t)-D_\nu(x)]|^pdx\nonumber\\
\leq
Kt^p\left\{\sum\limits_{m=1}^{n-1}m^{-2}\left(\sum\limits_{\nu=1}^m\nu^2|\Delta
a_\nu|\right)^p+\sum\limits_{m=1}^{n-1}m^{p-2}\left(\sum\limits_{\nu=m+1}^n\nu|\Delta
a_\nu|\right)^p\right\}.
\end{eqnarray*}
Thus, by Lemma 5 and Lemma 6, we obtain that
 \bll
I_{21}^p\leq
K\left(n^{-p}\sum\limits_{\nu=1}^{n-1}\nu^{2p-2}a_\nu^p+K_2
\sum\limits_{\nu=n}^{\infty}\nu^{p-2}a_\nu^p\right).\ell
In a way similar to the treatment of (3.9), we can easily deduce that
$$n^{p-1}a^{p}_{j(n+1)}\leq Kn^{-p}
\su_{\nu=1}^{n-1}\nu^{2p-2}a_\nu^p+\su_{\nu=n}^{\infty}\nu^{p-2}a_\nu^p,\;\;j=1,2,4,$$
with applying Lemma 4, we achieve that \bll
I_{22}&\leq&\left\{\int_{\pi/2n}^{\pi+\pi/2n}\left|\su_{\nu=n+1}^\infty|\Delta
a_\nu||D_\nu(x)|\right|^pdx\right\}^{1/p}\\
&\leq&K\left|\su_{\nu=n+1}^\infty|\Delta
a_\nu|\right|^p\left\{\int_{\pi/2n}^\infty
x^{-p}dx\right\}^{1/p}\\
&\leq&Kn^{1-1/p}\left(a_{n+1}+a_{2n+2}+a_{4n+4}+\su_{k=n+1}^\infty\frac{a_k}{k}\right)\\
&\leq&Kn^{-1}\left(\su_{\nu=1}^{n-1}\nu^{2p-2}a_\nu^p\right)^{1/p}+K
\left(\su_{\nu=n}^{\infty}\nu^{p-2}a_\nu^p\right)^{1/p}\\
&&+K\left(n^{p-1}\left(\su_{k=n+1}^{2n}\frac{a_k}{k}
\right)^p+n^{p-1}\left(\su_{k=2n}^\infty\frac{a_k}{k}\right)^p\right)^{1/p}\\
&\leq&Kn^{-1}\left(\su_{\nu=1}^{n-1}\nu^{2p-2}a_\nu^p\right)^{1/p}+K
\left(\su_{\nu=n}^{\infty}\nu^{p-2}a_\nu^p\right)^{1/p}\\
&&+K\left(n^{-1}\left(\su_{k=n}^{2n}a_k\right)^p+\su_{\nu=n+1}^{2n}\nu^{p-2}
\left(\su_{k=\nu}^\infty\frac{a_k}{k}
\right)^p\right)^{1/p}\\
&\leq&Kn^{-1}\left(\su_{\nu=1}^{n-1}\nu^{2p-2}a_\nu^p\right)^{1/p}+K
\left(\su_{\nu=n}^{\infty}\nu^{p-2}a_\nu^p\right)^{1/p}\\
&&+K\left(n^{p-2}\su_{k=n}^{2n}a_k^p+\su_{\nu=n+1}^\infty\nu^{p-2}
\left(\su_{k=\nu}^\infty\frac{a_k}{k} \right)^p\right)^{1/p}\\
&\leq&Kn^{-1}\left(\su_{\nu=1}^{n-1}\nu^{2p-2}a_\nu^p\right)^{1/p}+K
\left(\su_{\nu=n}^{\infty}\nu^{p-2}a_\nu^p\right)^{1/p}\\
&&+K\left(\su_{k=n}^{2n}k^{p-2}a_k^p+\su_{\nu=n+1}^\infty\nu^{p-2}
\left(\su_{k=\nu}^\infty\frac{a_k}{k} \right)^p\right)^{1/p}\\
&\leq&Kn^{-1}\left(\su_{\nu=1}^{n-1}\nu^{2p-2}a_\nu^p\right)^{1/p}+K
\left(\su_{\nu=n}^{\infty}\nu^{p-2}a_\nu^p\right)^{1/p}.\\
&&\hspace{6cm}\mbox{(by (4.1))}
 \ell
  Altogether, the above estimates for $I_1$ and $I_2$ complete
Theorem 1.

\vspace{3mm} {\bf Proof of Theorem 2.} By (3.8),
$$\sum\limits_{k=[n/2]}^{2n}a_k\geq K
na_n,\hspace{.4in}(4.2)$$
and combining (4.2) with (3.1) of Lemma 1, we have (2.4).

By estimate (3.2) of Lemma 1, (2.7) will be proved if the
following inequality
\begin{eqnarray*}
\sum\limits_{n=1}^\infty\lambda(n)\left(\sum\limits_{k=n}^\infty|\Delta
a_k|\right)^r\leq K\sum\limits_{n=1}^\infty\lambda(n)a_n^r
\end{eqnarray*}
holds. Furthermore, with the help of Lemma 4, what we really need to establish is that
\begin{eqnarray*}
\sum\limits_{n=1}^\infty\lambda(n)\left(\sum\limits_{k=n}^\infty\frac{a_k}{k}\right)^r\leq
K\sum\limits_{n=1}^\infty\lambda(n)a_n^r.
\end{eqnarray*}
In fact, if $r=1$, by exchanging the order of summation and applying (2.5),
we have
\begin{eqnarray*}\sum\limits_{n=1}^\infty\lambda(n)
\sum\limits_{k=n}^\infty\frac{a_k}{k}&=&\sum\limits_{n=1}^\infty\frac{a_n}{n}\sum\limits_{k=1}^n
\lambda(k)k^{1/p-1}k^{1-1/p} \\
&\leq&
\sum\limits_{n=1}^\infty\frac{a_n}{n}n^{1-1/p}\sum\limits_{k=1}^n
\lambda(k)k^{1/p-1}\leq
K\sum\limits_{n=1}^\infty\lambda(n)a_n;\end{eqnarray*}
 if $r>1$, by H\"{o}lder's inequality and (2.5), we still have
\begin{eqnarray*}
\sum\limits_{n=1}^\infty\lambda(n)
\left(\sum\limits_{k=n}^\infty\frac{a_k}{k}\right)^r&\leq&\sum\limits_{n=1}^\infty\lambda(n)
\left(\sum\limits_{k=n}^\infty\frac{1}{k^{
1-\frac{(1-p)r}{p(r-1)}}}
\right)^{r-1}\sum\limits_{k=n}^\infty\frac{a_k^r}{k^{1+\frac{(1-p)r}{p}}}\\
&\leq&K\sum\limits_{k=1}^\infty\frac{a_k^r}{k^{1+\frac{(1-p)r}{p}}}\sum\limits_{n=1}^k\lambda(n)n^{\frac{r}{p}-r}
\leq \sum\limits_{n=1}^\infty\lambda(n)a_n^r.
\end{eqnarray*}

{\bf Proof of Theorem 3.} Most of the proof can be proceeded as the
corresponding part of Leindler [7] word by word, we omit the details
here.

\vspace{3mm} {\bf Proof of Theorem 4.} As Leindler [7] pointed
out, what we need to do is to verify that $(2.8)\Rightarrow$
(i) and that (ii) $\Rightarrow (2.8)$. By applying Abel's tranformation
$$\su_{\nu=n}^\infty\lambda_\nu u_\nu=\su_{\nu=n}^\infty(\lambda_\nu-\lambda_{\nu+1})
\su_{k=1}^\nu u_k-\lambda_n\su_{k=1}^{n-1}u_k$$ with
$\lambda_\nu=\nu^{-p}$ and $u_\nu=\nu^{2p-2}a_\nu^p$,  and also by (2.5), we evidently have
$$\su_{\nu=n}^\infty\nu^{p-2}a_\nu^p=\su_{\nu=n}^\infty\nu^{-p}\nu^{2p-2}a_\nu^p
\leq Kn^{-p},$$ thus,  the second term in (2.3) is
not larger than $Kn^{-1}$. Altogether, by Theorem 1, it means that $f\in\Lambda_p$.

Let $f(x)$ be the sum function of, say, the series $\sum_{n=1}^{\infty}\limits a_n\sin nx$, and
set
$$F(x):=\int_{0}^xf(t)dt=\sum\limits_{n=1}^\infty n^{-1}a_n(1-\cos nx).$$
An standard argument  yields that
$$F(\pi/(2n))=2\sum\limits_{k=1}^\infty\frac{a_k}{k}
\sin^2\frac{k\pi}{4n}\geq
\frac{K}{n}\sum\limits_{k=[n/2]}^{2n}a_k,$$ so that $F(\pi/(2n))\geq
Ka_n$ by (4.2). Set
$$G(x):=\int_0^xdt\int_0^t|f'(u)|du.$$
Obviously, $F(x)\leq G(x)$. Hence, applying Hardy's
inequality ([12]) twice, we obtain that \bll\su_{n=2}^\infty
n^{2p-2}a_n^p&\leq&\su_{n=2}^\infty n^{2p-2}G^{p}(\pi/(2n))\leq
K\su_{n=2}^\infty n^{2p-2}G^{p}(\pi/n)\\
&\leq&K\su_{n=2}^\infty\int_{\pi/n}^{\pi/(n-1)}\left[\frac{G(x)}{x}\right]^{p}x^{-p}dx
\leq K\int_0^\pi\left[\frac{G(x)}{x}\right]^{p}x^{-p}dx\\
&\leq&K\int_0^\pi\left(\int_0^x|f'(t)|dt\right)^px^{-p}dx
\leq K\int_0^\pi|f'(x)|^pdx<\infty.\ell

 {\bf Proof
of Theorem 5.} Set $T_{m,2n}(x):=\su_{\nu=m}^{2n}\cos\nu x$, then
([1] or [7]) \bll
I_{m,2n,t}:=\int_{-\pi}^\pi\left(2f(x)-f(x+t)-f(x-t)\right)T_{m,2n}(x)dx=4\pi\su_{\nu=m}^{2n}
a_\nu\sin^2\frac{1}{2}\nu t.\hspace{.4in}(4.3)\ell Taking $t=\pi/n$ and $m=[n/2]$ in
(4.3), then applying (4.2) again, we have \bll\su_{\nu=m}^{2n}
a_\nu\sin^2\frac{1}{2}\nu t\geq K\su_{\nu=m}^{2n} a_\nu\geq
na_n. \hspace{.4in}(4.4)\ell On the other hand, we have \bll
\int_{-\pi}^\pi|T_{m,2n}(x)|^qdx\leq
K\left\{\int_0^{\pi/(2n)}n^qdx+\int_{\pi/(2n)}^\pi
x^{-q}dx\right\}\leq Kn^{q-1}.\ell By H\"{o}lder's inequality,
it follows that \bll I_{m,2n,\pi/n}\leq
Kn^{1/p}\left\{\int_{-\pi}^\pi\left|f\left(x+\frac{\pi}{n}\right)+
f\left(x-\frac{\pi}{n}\right)-2f(x)\right|^pdx\right\}^{1/p}\leq
Kn^{1/p}\omega_p^*(f,\pi/n).\ell A combination of (4.3) and (4.4) leads to \bll\omega_p^*(f,1/n)\geq Kn^{1-1/p}a_n.\hspace{.4in}(4.5)\ell Therefore, from that
$\omega_p^*(f,1/n)\leq Kn^{-1}$ and by (4.5), we get $a_n\leq
Kn^{-2+1/p}$, whence by Theorem 1, it follows that
$$\omega_p(f,1/n)\leq Kn^{-1}(\log n)^{1/p},$$
and (2.9) is done.

\begin{center}
{\Large\bf  References}
\end{center}
\begin{enumerate}

\item \bf S. Aljan$\breve{c}$i\'{c}, \it On the integral moduli of
continuity in $L^p (1<p<\infty)$ of Fourier series with monotone
coefficients, \rm Proc. Amer. Math. Soc., 17(1966), 287-294.

\item \bf R. J. Le and S. P. Zhou , \it A new condition for the
uniform convergence of certain trigonometric series, \rm Acta
Math. Hungar., 108(2005), 161-169.

 \item \bf L. Leindler, \it Generalization of
inequalities of Hardy and Littlewood, \rm Acta Sci. Math.
(Szeged), 31(1970), 279-285.

\item \bf L. Leindler, \it On cosine series with positive
coefficients, \rm Acta Math. Hungar., 22(1971),
397-406.

\item \bf L. Leindler, \it On the uniform convergence and
boundedness of a certain class of sine series, \rm Anal. Math.,
27(2001), 279-285.

\item \bf L. Leindler, \it A new class of numerical sequences and
its applications to sine and cosine series, \rm Anal. Math.,
28(2002), 279-286.

\item \bf L. Leindler, \it Relations among Fourier series and
sum-functions, \rm  Acta Math. Hungar., 104(2004), 171-183.

\item\bf D. S. Yu and S. P. Zhou, \it A generalization of
monotonicity and applications, \rm submitted to Acta Math.
Hungar..

\item\bf D. S. Yu and S. P. Zhou, \it  On belonging of
trigonometric series to Ba spaces, \rm to appear.

\item \bf S. P. Zhou and R. J. Le, \it A remark on ``two-sided"
monotonicity condition: an application to $L^p$ convergence, \rm
Acta Math. Hungar. accepted.

\item \bf S. P. Zhou and R. J. Le, \it A new condition for
the uniform convergence in Fourier analysis, \rm Adv. Math. (Beijing), 33(2004), 567-569.

\item\bf A. Zygmund, \it Trigonometric Series, 2nd. Ed.,
Vol.I, \rm Cambridge Univ. Press, Cambridge, 1959.

\end{enumerate}

\vspace{4mm}
\begin{flushleft}

\bf Yu Dansheng  \\
\rm Institute of Mathematics\\
Zhejiang Sci-Tech University \\
Xiasha Economic Development Area\\
Hangzhou Zhejiang 310018 China\\
e-mail: danshengyu@yahoo.com.cn

\vspace{3mm}

\bf Zhou Songping \\
\rm Institute of Mathematics\\
Zhejiang Sci-Tech University \\
Xiasha Economic Development Area\\
Hangzhou Zhejiang 310018 China\\
e-mail: szhou@zjip.com
\end{flushleft}

\end{document}